\documentclass[10pt]{article}
\usepackage{latexsym,amsmath,amssymb,graphics,amscd}

\makeatletter
\@addtoreset{figure}{section}
\def\thefigure{\thesection.\@arabic\c@figure}
\def\fps@figure{h,t}
\@addtoreset{table}{bsection}

\def\thetable{\thesection.\@arabic\c@table}
\def\fps@table{h, t}
\@addtoreset{equation}{section}

\makeatother

\begin{document}

\newtheorem{theorem}{Theorem}[section]
\newtheorem{definition}[theorem]{Definition}
\newtheorem{lemma}[theorem]{Lemma}
\newtheorem{remark}[theorem]{Remark}
\newtheorem{proposition}[theorem]{Proposition}
\newtheorem{corollary}[theorem]{Corollary}
\newtheorem{example}[theorem]{Example}
\newtheorem{examples}[theorem]{Examples}

\newcommand{\bfi}{\bfseries\itshape}

\newsavebox{\savepar}
\newenvironment{boxit}{\begin{lrbox}{\savepar}
\begin{minipage}[b]{15.8cm}}{\end{minipage}
\end{lrbox}\fbox{\usebox{\savepar}}}

\makeatletter
\title{{\bf Poisson reduction}}
\author{Juan-Pablo Ortega$^{1}$ and  Tudor S. Ratiu$^{2}$}
\addtocounter{footnote}{1}
\footnotetext{Centre National de la Recherche Scientifique,
D\'epartement de Math\'ematiques de Besan\c con,
Universit\'e de Franche-Comt\'e.
UFR des Sciences et Techniques.
16, route de Gray.
F-25030 Besan\c con cedex. France. {\texttt
Juan-Pablo.Ortega@math.univ-fcomte.fr}. }
\addtocounter{footnote}{1}
\footnotetext{Section de Math\'ematiques, 
\'Ecole Polytechnique F\'ed\'erale de Lausanne,  CH-1015 Lausanne,
Switzerland. {\texttt Tudor.Ratiu@epfl.ch}.}
\date{}
\makeatother
\maketitle

\begin{abstract}
This encyclopedia article briefly reviews without proofs some of the main results
in Poisson reduction. The article recalls most the necessary prerequisites to
understand the main results.
\end{abstract}

The Poisson reduction techniques allow the construction of new Poisson 
structures out of a given one by combination of two operations: 
{\bfi  restriction} to submanifolds that satisfy certain compatibility 
assumptions and  passage to a {\bfi  quotient space} where
certain degeneracies have been eliminated. For certain kinds of reduction
it is necessary to pass first to a submanifold and then take a quotient.
Before making this more explicit we introduce the notations that will be used in this article. All manifolds
in this article are finite dimensional.

\medskip

\noindent  {\bf Poisson manifolds.} A {\bfi Poisson
manifold\/} is a pair $(M,\,\{\cdot,\cdot\})$, where $M$ is a manifold 
and $\{\cdot,\cdot\}$ is a bilinear operation on $C^{\infty}(M)$
such that $(C^{\infty}(M),\,\{\cdot,\cdot\})$ is a Lie algebra
and $\{\cdot,\cdot\}$ is a derivation (that is, the Leibniz
identity holds) in each argument. The
pair $(C^{\infty}(M),\,\{\cdot,\cdot\})$ is also called a {\bfi Poisson algebra}.
The functions in the center $\mathcal{C}(M)$
of the Lie algebra $(C^{\infty}(M),\,\{\cdot ,\cdot \})$ are
called {\bfi Casimir functions\/}.
From the natural isomorphism between derivations on $C^{\infty}(M)$
and  vector fields on $M$, it follows that
each $h\in C^{\infty}(M)$  induces a vector field on $M$ via the expression
$X_h = \{\cdot , h\}$,
called the {\bfi Hamiltonian vector field\/}
\index{Hamiltonian!vector field}
\index{vector!field!Hamiltonian}%
associated to  the {\bfi Hamiltonian function\/}
\index{Hamiltonian!function}
\index{function!Hamiltonian}%
$h$. 
The triplet
$(M,\,\{\cdot,\cdot \},\,h)$ is called a {\bfi Poisson dynamical 
system\/}.
Any Hamiltonian system on a symplectic manifold is a
Poisson dynamical system relative to the Poisson bracket
induced by the symplectic structure.
Given a Poisson dynamical system $(M,\,\{\cdot,\cdot \},\,h)$, 
its {\bfi integrals of motion\/} 
\index{integral!of the motion}
or {\bfi conserved
quantities\/}
\index{conserved quantity}
are defined as the centralizer of $h$ in $(C^{\infty}(M),\,\{\cdot
,\cdot \})$ that is, the subalgebra of $(C^{\infty}(M),\,\{\cdot
,\cdot \})$ consisting of the functions $f\in C^{\infty}(M)$ such
that $\{f,\,h\}=0$. Note that the terminology is justified
since, by Hamilton's equations in Poisson bracket form, we have 
$\dot{f}=X_h [f]=\{f,\,h\}=0$,
that is, $f$ is constant on the flow of $X_h$. A smooth mapping $\varphi:M_1\rightarrow M_2$, between the two 
Poisson manifolds $(M_1,\,\{\cdot,\cdot\}_1)$ and
$(M_2,\,\{\cdot,\cdot\}_2)$ is called {\bfi canonical\/} 
\index{canonical!map}
\index{map!canonical}%
or
{\bfi Poisson\/} 
if for all $g,\,h\in C^{\infty}(M_2)$ we have
$
\varphi^\ast \{g,\, h\}_2=\{\varphi^\ast  g,\,\varphi^\ast 
g\}_1\,$.
If $\varphi:M_1\rightarrow M_2$ is a smooth map between two 
Poisson manifolds $(M_1,\,\{\cdot,\cdot\}_1)$ and
$(M_2,\,\{\cdot,\cdot\}_2)$ then
$\varphi$ is a Poisson map if and only if $T\varphi\circ 
X_{h\circ \varphi}=X_h\circ \varphi$ for any $h \in
C^\infty(M_2)$, where $T \varphi: TM_1 \rightarrow TM_2 $ denotes the
tangent map (or derivative) of $\varphi$.

Let $(S,\{ \cdot , \cdot \}^S)$
and $(M,\{ \cdot , \cdot \} ^M)$ be two Poisson manifolds such that $S \subset M $ and
the inclusion $i _S:S\hookrightarrow M $  is an immersion. The Poisson manifold $(S,\{ \cdot , \cdot \}^S)$
is called a {\bfi  Poisson submanifold} 
of $(M,\{ \cdot , \cdot \} ^M)$ if  $i _S $ is a canonical map. 
An immersed submanifold $Q $ of $M$  is called a {\bfi  quasi Poisson
submanifold} of 
$(M,\{ \cdot , \cdot \} ^M)$ if for any $q \in Q $, any open neighborhood $U$ of $q$
in $M$, and any $f \in C^{\infty} (U) $ we have 
$
X _f(i _Q (q)) \in T _qi _Q(T _qQ)$,
where $i _Q:Q\hookrightarrow M $ is the inclusion and $X _f $ is the Hamiltonian vector
field of $f$ on $U$ with respect to the  Poisson bracket of $M$ restricted to $U$.
If $(S,\{ \cdot , \cdot \}^S)$ is a Poisson submanifold of $(M,\{ \cdot , \cdot
\} ^M)$ then there is no other bracket $\{ \cdot , \cdot \}' $ on $S$ making the
inclusion $i:S
\hookrightarrow M $ into a canonical map. 
If $Q$ is a quasi Poisson
submanifold of $(M,\{ \cdot , \cdot \})$ then there exists a unique Poisson
structure
$\{
\cdot , \cdot \} ^Q $ on $Q$ that makes it into a Poisson submanifold of $(M,\{ \cdot ,
\cdot \})$ but this Poisson structure may be different from the given one on $Q$. Any Poisson
submanifold is quasi Poisson but the converse is not  true in general.

\medskip
 
\noindent  {\bf The Poisson tensor and symplectic leaves.} The derivation property of the
Poisson bracket implies that for  any two functions $f,\,g\in C^{\infty}(M)$, the value
of the bracket $\{f,\,g\}(z)$ at an arbitrary point $z\in M$ (and
therefore $X_f(z)$ as well), depends on $f$ only through $\mathbf{d}
f(z)$  which allows us to define  a contravariant
antisymmetric two-tensor $B\in\Lambda^2(T^\ast  M)$, called the {\bfi
Poisson tensor\/},  by 
$B(z)(\alpha_z,\,\beta_z)=\{f,\,g\}(z)$,
where $\mathbf{d} f(z)=\alpha_z \in T^\ast _z M$ and $\mathbf{d} g(z)=\beta_z\in T^\ast _z M$. 
The vector bundle map $B^\sharp:T^\ast  M\rightarrow TM$ over the identity
naturally associated to $B$  is defined by
$B(z)(\alpha_z,\,\beta_z)=
\langle\alpha_z,\,B^\sharp(\beta_z)\rangle$.
Its range $D:=B^\sharp(T^\ast  M) \subset TM$ is called the
{\bfi characteristic distribution} of  $(M, \{ \cdot , \cdot \}) $ since $D$ is  a
generalized smooth integrable distribution. Its maximal integral leaves
are called the {\bfi  symplectic leaves} of $M$ for they carry a symplectic
structure that makes them into Poisson submanifolds. As integral leaves of an integrable distribution, the
symplectic leaves $\mathcal{L}  $ are {\bfi  initial submanifolds } of  $M$, that
is, the inclusion $i:\mathcal{L}\hookrightarrow M $ is an injective immersion
such that  for any smooth manifold $P$, an arbitrary map $g :P \rightarrow 
\mathcal{L}  $  is smooth if  and only if $i \circ g: P  \rightarrow  M  $ 
is smooth.

\section{Poisson reduction}

\noindent {\bf Canonical Lie group actions.} Let $(M,\,\{\cdot,\cdot\})$
be a Poisson manifold and let $G$ be a Lie group acting
canonically on $M$ via the map $\Phi:G \times  M \rightarrow  M $. An action is
called {\bfi  canonical} if  for any
$h
\in G
$ and $f,g
\in  C^\infty(M)$ one has  
\[
\{f \circ \Phi _h, g \circ  \Phi _h\}=\{f,g\} \circ \Phi_h.
\]
If the
$G$-action is free and proper then the
orbit space $M/G$ is a smooth regular quotient manifold. Moreover, it  is also a Poisson
manifold with  the Poisson bracket $\{\cdot,\cdot\}^{M/G}$, uniquely characterized by
the relation
\begin{equation}
\label{characterization reduced free bracket}
\{f,\,g\}^{M/G}(\pi (m))=\{f\circ\pi,\,g\circ\pi\}(m),
\end{equation}
for any $m \in M $ and where  $f,\,g:M/G\rightarrow\mathbb{R}$ are two arbitrary
smooth  functions. This bracket is appropriate for the reduction of Hamiltonian dynamics
in the sense that if
$h \in C^\infty(M)^{G}$ is  a $G$-invariant smooth function on
$M$, then the Hamiltonian flow $F_t$ of $X_h$ commutes with the
$G$-action, so it induces a flow $F_t^{M/G}$ on $M/G$
that  is Hamiltonian on
$(M/G,\,\{\cdot,\cdot\}^{M/G})$ for the {\bfi reduced Hamiltonian}
function $[h]\in C^{\infty}(M/G)$ defined by
$
[h]\circ\pi=h$.

If the Poisson manifold $(M,\,\{\cdot,\cdot\})$ is actually symplectic with
form $\omega$ and the $G$-action has an associated momentum map $\mathbf{J}:
M \rightarrow \mathfrak{g}^\ast $, then the symplectic leaves of $(M/G, \{\cdot ,
\cdot\}^{M/G})$ are given by the spaces $ \left( M _{\mathcal{O}_{\mu}} ^c:=G \cdot
\mathbf{J}^{-1}(\mu) ^c/G , \omega _{\mathcal{O}_{\mu}}^c\right)$, where
$\mathbf{J}^{-1}(\mu) ^c  $ is a connected component of the fiber
$\mathbf{J}^{-1}(\mu) $ and $\omega _{\mathcal{O}_{\mu}} ^c $ is the restriction to
$M _{\mathcal{O}_{\mu}}^c  $ of the symplectic form $ \omega _{\mathcal{O}_{\mu}}$ of
the symplectic orbit reduced space $M _{\mathcal{O}_{\mu}} $  (see~\cite{symplectic
reduction encyclopedia}). If, additionally,
$G$ is compact, $M$  is connected, and the momentum map $\mathbf{J}$ is proper then $M
_{\mathcal{O}_{\mu}} ^c= M _{\mathcal{O}_{\mu}} $.

\medskip

In the remainder of this section we characterize the situations in which new Poisson
manifolds can be obtained out of a given one by a combination of restriction to a
submanifold and passage to the quotient with respect to an equivalence relation that
encodes the symmetries of the bracket. 

\begin{definition}
\label{poisson distribution weak}
Let $(M,\,\{\cdot ,\cdot\})$ be a Poisson manifold and 
$D\subset TM$  a smooth distribution on $M$. The
distribution $D$ is called {\bfi Poisson\/} or {\bfi canonical\/},
\index{Poisson!distribution}%
\index{canonical!distribution}%
\index{distribution!Poisson}%
\index{distribution!canonical}%
if the condition $\mathbf{d} f|_D=\mathbf{d} g|_D=0$, for any $f,\,g\in
C^{\infty}(U)$ and any open subset $U \subset P $, 
implies that $\mathbf{d}\{f,\,g\}|_D=0$
\end{definition}

Unless strong regularity assumptions are invoked, the passage to the leaf
space of a canonical distribution destroys the smoothness of the quotient
topological space. In such situations the Poisson algebra of functions 
is too small and the notion of  {\bfi  presheaf of Poisson algebras} is needed.
\begin{definition}
\label{presheaf of Poisson algebras definition}
Let $M$ be a topological space with a presheaf ${\cal F}$ of smooth functions. A {\bfi 
presheaf of Poisson algebras} 
on $(M, \mathcal{F}) $ is a map $\{ \cdot , \cdot \} $ that
assigns to each open set $U \subset M $ a bilinear operation $\{ \cdot , \cdot \} _U :
\mathcal{F} (U) \times \mathcal{F}(U) \rightarrow \mathcal{F}(U)$ such that the pair
$(\mathcal{F}(U),\{ \cdot , \cdot \} _U)$ is a  Poisson algebra. A presheaf of Poisson
algebras is denoted as a triple $(M, \mathcal{F}, \{ \cdot , \cdot \})$.
The presheaf of Poisson algebras $(M, \mathcal{F}, \{ \cdot , \cdot \})$ is
said to be {\bfi  non-degenerate}
if the following condition holds: if $f \in \mathcal{F}(U)$ is such that $\{f,g\}_{U
\cap V}=0 $, for any $g \in \mathcal{F}(V) $ and any open set of $V$, then $f $ is
constant on the connected components of $U$.
\end{definition}
Any  Poisson manifold $(M, \{ \cdot , \cdot \})$ has a natural presheaf of Poisson
algebras on its presheaf  of smooth functions that associates to any open subset $U $ of
$M$ the restriction $\{ \cdot , \cdot \}| _U $ of $\{ \cdot , \cdot \} $ to $C^{\infty}
(U)\times C^{\infty} (U)$.

\begin{definition}
\label{stratified set}
Let $P$ be a topological space and ${\cal Z} = \{S_i\}_{i \in I}$ a locally
finite partition of $P$ into smooth manifolds $S _i \subset P $, $i \in
I$, that are locally closed topological subspaces of $P$ (hence
their manifold topology is the relative one induced by $P$). The
pair $(P, {\cal Z})$ is called a {\bfi  decomposition} of $P$ with
{\bfi  pieces} in ${\cal Z}$, or a {\bfi decomposed space}, if the
following {\bfi frontier condition} holds:
\begin{description}
\item [(DS)] If $R, S \in {\cal Z}$ are such that $R \cap \bar{S} 
\neq \emptyset $, then $R \subset \bar{S} $. In this case we write
$R \preceq S $. If, in addition, $R \neq S$ we say that
$R$ is {\bfi  incident } to $S$ or that it is a {\bfi  boundary
piece} of $S$ and write $R \prec S$.  
\end{description}
\end{definition}

\begin{definition}
\label{smooth adapted distribution decomposition}
Let $M$ be a differentiable manifold and $S\subset M$  a 
decomposed subset of $M$. Let 
$\{S_i\}_{i\in I}$ be the pieces of this decomposition. The topology of $S$ is not
necessarily the relative topology as a subset of $M$.  Then $D\subset TM|_S$
is called a {\bfi smooth distribution on $S$ adapted to the decomposition\/} 
\index{distribution!adapted to a decomposition}%
\index{decomposition!distribution adapted to a}%
\index{distribution!smooth}%
\index{smooth!distribution}%
$\{S_i\}_{i\in I}$, if
$D\cap TS_i$ is a smooth distribution on $S_i$ for all $i\in I$. The
distribution $D$ is said to be {\bfi integrable\/} 
\index{distribution!integrable}%
\index{integrable!distribution}%
if $D\cap TS_i$
is integrable for each $i\in I$.
\end{definition}
In the situation described by the previous definition and if $D$ is integrable,  the 
integrability of the distributions $D_{S _i}:=D\cap TS_i$ on $S_i$ allows
us to partition each $S_i$ into the corresponding maximal integral manifolds. 
Thus, there is an equivalence relation  on $S_i$ whose 
equivalence classes are precisely these maximal integral
manifolds. Doing this on each $S_i$, we obtain an equivalence 
relation $D_S$ on the whole set $S$ by taking the union of the
different  equivalence classes corresponding to all the $D_{S _i}$.
Define the quotient space $S/D _S$ by
\[S/D _S:=\bigcup_{i\in I} S_i/D_{S _i}\]  
and let $\pi_{D_S }: S \rightarrow  S/ D_S $ be the natural projection.

\medskip

\noindent {\bf The presheaf of smooth functions on $S/D_S$.} Define the
presheaf of smooth functions
$C^{\infty} _{S/ D_S} $ on $S/ D _S $ as the map that associates to any open subset $V $
of $S/ D _S$ the set of functions $C^{\infty} _{S/ D_S} (V)$ characterized
by the following property:
$f 
\in C^{\infty} _{S/ D_S} (V)$ if and only if for any $z \in V $  there
exists $m
\in \pi _{D_S} ^{-1} (V)
$, $U _m $ open neighborhood of $m$ in  $M$, and
$F \in C^{\infty} (U _m) $ such that 
\begin{equation}
\label{one more presheaf}
f \circ  \pi_{D _S}|_{\pi _{D_S} ^{-1} (V)\cap U _m}=F|_{\pi _{D_S} ^{-1} 
(V)\cap U_m};
\end{equation}
$F$ is called a {\bfi  local extension}
\index{extension!local}%
\index{local!extension}%
of $f \circ \pi_{D_S} $ at the point $m \in \pi_{D_S} ^{-1} (V) $. When the
distribution $D$ is trivial, the presheaf $C^{\infty} _{S/ D_S}$ coincides with
the presheaf of {\bfi  Whitney smooth functions} $C^{\infty}_{S,M}$ on $S$
induced by the smooth functions on  $M$.

The presheaf $C^{\infty}_{S/ D_S}$ is said to have the $(D, D_S)$-{\bfi 
local extension property}
\index{local!extension!property}%
\index{extension!property!local}%
when the topology of $S$ is stronger than the relative topology and, at the same time,
the local extensions of
$f  \circ \pi_{D_S} $ defined in~(\ref{one more presheaf}) can always be chosen to
satisfy 
\[
\mathbf{d} F (n) |_{D(n)}=0, \quad \mbox{for any}\quad n \in  
\pi _{D_S} ^{-1} (V)\cap U_m;
\] 
$F$ is called a {\bfi  local $D$-invariant extension}
of $f \circ \pi_{D_S} $ at the point $m \in \pi_{D_S} ^{-1} (V) $. 
If $S$ is a smooth embedded submanifold of $M$ and 
$D_S $ is a smooth, integrable, and regular distribution on $S$, then the
presheaf
$C^{\infty}_{S/ D_S}$ coincides with the presheaf of smooth functions on $S / D_S $ 
when considered as a regular quotient manifold.

The following definition spells out what we mean by obtaining a bracket via reduction.

\begin{definition}
Let  $(M, \{ \cdot , \cdot \})$ be a Poisson manifold, $S$  a decomposed subset of $M$,
and $D \subset  TM| _S $  a Poisson integrable generalized distribution adapted to the
decomposition of $S$. Assume that  $C^{\infty} _{S/ D_S}$ has the $(D,
D_S)$-local extension property. Then $(M, \{ \cdot , \cdot \}, D, S)$
is said to be {\bfi  Poisson reducible} 
\index{Poisson!reducible}%
\index{reducible!Poisson}%
if  $(S/ D_S, C^{\infty}_{S/ D_S},\{ \cdot , \cdot \}^{S/ D_S})$ is a well
defined presheaf of Poisson algebras where, for any open set $V \subset S/ D_S$, the
bracket $\{
\cdot , \cdot \}_V ^{S/ D_S}: C^{\infty} _{S/ D_S} (V) \times C^{\infty} _{S/ D_S} (V)
\rightarrow C^{\infty} _{S/ D_S} (V) $ is given by
\[
\{
f , g \}_V ^{S/ D_S}(\pi_{D_S} (m)):=\{F, G\} (m),
\]
for any $m \in \pi_{D_S}^{-1} (V)$ for local $D$-invariant extensions 
$F,G $ at $m$ of $f \circ \pi_{D_S} $ and $g \circ
\pi_{D_S} $, respectively.
\end{definition}

\begin{theorem}
\label{reduction distribution singular only one}
Let $(M, \{ \cdot , \cdot \}) $ be a Poisson manifold with associated Poisson tensor $B
\in  \Lambda ^2(T ^\ast M)$, $S$  a decomposed space, and $D \subset  TM| _S
$  a Poisson integrable generalized distribution adapted to the decomposition
of $S$ (see Definitions~{\rm \ref{smooth adapted distribution
decomposition}} and~{\rm \ref{poisson distribution weak}}). Assume that 
$C^{\infty} _{S/ D_S}$ has the $(D, D_S)$-local extension property. Then
$(M, \{
\cdot , \cdot \}, D, S)$ is Poisson reducible if for any
$m \in S  $
\begin{equation}
\label{Poisson reducibility distribution singular}
B^{\sharp}(\Delta _m) \subset  \left[ \Delta _m ^S  \right] ^{\circ}
\end{equation}
where $\Delta_m:=\{ \mathbf{d} F (m)\mid F \in C^{\infty} (U _m), \mathbf{d}F
(z)|_{D (z)}=0, \text{ for all }z \in U _m\cap S,\text{ and for any open neighborhood}$ 
$\text{$U _m  $ of $m$ in $M$}\}
$ and $\Delta _m^S:=\{ \mathbf{d} F (m) \in \Delta _m\mid F|_{U _m\cap V _m}\text{ is
constant for an open neighborhood }$  $U _m    \text{ of }   m  \text{ in }$  $M \text{
and an open neighborhood
$V _m
$ of $m$ in
$S$}\}
$.
\end{theorem} 
If $S$ is endowed with the relative topology then $\Delta _m^S:=\{ \mathbf{d} F (m) 
\in \Delta _m\mid F|_{U _m\cap V _m}\text{ is
constant}$ $\text{for an open}$ $\text{neighborhood }$  $U _m    \text{ of }   m  \text{
in }  M\}
$.

\medskip

\noindent {\bf Reduction by regular canonical distributions.} Let $(M,\{ \cdot , \cdot
\})$ be  a Poisson manifold and $S$   an embedded submanifold of $M$. Let $D\subset TM|
_S $  be a subbundle of the tangent bundle of $M$ restricted to $S$ such that 
 $D_S:=D\cap TS  $  is a smooth, integrable, regular distribution  on
$S$ and $D$ is canonical.

\begin{theorem}
\label{reduction theorem by distributions regular case} 
With the above hypotheses, $(M, \{ \cdot , \cdot \}, D,
S) $ is Poisson reducible if and only if 
\begin{equation}
\label{classical reduction condition Marsden Ratiu}
B ^{\sharp}(D^{\circ})\subset TS + D.
\end{equation}
\end{theorem}

\section{Applications of the Poisson Reduction Theorem}

\noindent  {\bf Reduction of coisotropic submanifolds.} 
Let $(M,\{ \cdot , \cdot \})$ be a Poisson manifold with associated Poisson tensor $B
\in  \Lambda ^2(T ^\ast M)$ and $S$   an immersed smooth submanifold of  $M$.
Denote by $(TS)^\circ : =  \{ \alpha_s \in T ^\ast _s M \mid \langle
\alpha_s, v_s \rangle =  0, \text{~for~all~} s \in S, v_s \in T_s S\} \subset T
^\ast M$ the {\bfi conormal bundle\/}
\index{conormal bundle}%
\index{bundle!conormal}%
of the manifold $S $; it is a vector subbundle of $T ^\ast M|_S $.  The manifold $S $
is called {\bfi coisotropic\/}
\index{coisotropic!submanifold}%
\index{submanifold!coisotropic}%
if $B ^{\sharp}\left( (TS)^\circ \right) \subset TS $. 
 In the physics literature, coisotropic submanifolds appear
sometimes  under the name of {\bfi  first class constraints}.
The
following are equivalent:
\begin{description}
\item[(i)] $S$ is coisotropic;
\item [(ii)] if $f \in  C^\infty(M) $ satisfies  $f| _S\equiv 0 $ then $X _f| _S \in
\mathfrak{X} (S) $;
\item[(iii)] for any $s \in S$, any open neighborhood $U_s $ of $s $ in $M $, and any
function $g \in C^\infty(U_s)$ such that $X_g (s) \in T_s S $, if $f \in
C^\infty(U_s) $ satisfies $\{f, g\}(s) = 0 $, it follows that $X_f(s) \in T_sS $;
\item[(iv)] the subalgebra $\{f \in C^\infty(M) \mid f|_S \equiv 0 \}$ is a Poisson
subalgebra of $ (C^\infty(M), \{\cdot, \cdot \})$. 
\end{description}
The following proposition shows how to endow the coisotropic submanifolds of a Poisson
manifold with a Poisson structure by using the Reduction Theorem~\ref{reduction
distribution singular only one}.

\begin{proposition}
\label{everything works with coisotropic}
Let $(M,\{ \cdot , \cdot \})$ be a Poisson manifold with associated Poisson tensor $B
\in \Lambda ^2(T ^\ast M) $. Let $S$ be an embedded coisotropic submanifold of $M$ and
$D :=B ^{\sharp} ((TS)^{\circ})$. Then
\begin{description}
\item [(i)] $D= D\cap TS = D_S$ is a smooth generalized distribution on $S$.
\item  [(ii)] $D$ is integrable.
\item [(iii)] If $C^{\infty}_{S/D_S} $ has the $(D, D_S)$-local extension
property then $(M, \{ \cdot , \cdot \}, D, S)$ is  Poisson reducible. 
\end{description}
\end{proposition}

Coisotropic submanifolds usually appear as the level sets of integrals in 
involution. Let $(M,\{\cdot , \cdot \})$ be a Poisson manifold with Poisson
tensor $B$ and let $f_1, \ldots, f _k
\in  C^\infty(M) $ be $k$ smooth functions in {\bfi  involution},
\index{involutive!functions}%
that is,
$\{f _i, f _j\}=0$, for any $i,j \in \{1, \ldots, k\}$.
Assume that ${\bf 0} \in \mathbb{R}^k $ is a regular value of the function $F:=(f _1,
\ldots, f _k):M \rightarrow  \mathbb{R} ^k $ and let $S:= F ^{-1} (0) $. Since for
any
$s \in S $, ${\rm span}\{ \mathbf{d} f _1 (s), \ldots, \mathbf{d} f _k (s)\} \subset
\left( T _s S \right)^{\circ} $ and the dimensions of both sides of this inclusion are
equal it follows that
${\rm span}\{ \mathbf{d} f _1 (s), \ldots, \mathbf{d} f _k (s)\}=
\left( T _s S \right)^{\circ}$.
Hence
$B^{\sharp}(s)((T_{s}S)^{\circ})= {\rm span} \left\{ X_{f _1}(s),  \ldots,
X_{f _k}(s)\right\}$ and  $B^{\sharp}(s)((T_{s}S)^{\circ}) \subset T_s S$ 
by the involutivity of the components of $F$. Consequently, $S$ is a
coisotropic submanifold of $(M,\{ \cdot ,
\cdot \})$.

\medskip

\noindent {\bf Cosymplectic submanifolds and Dirac's constraints formula.}
The Poisson Reduction Theorem~\ref{reduction theorem
by distributions regular case} allows us to define Poisson structures on certain
embedded submanifolds that are not Poisson submanifolds. 

\begin{definition}
Let $(M, \{ \cdot , \cdot \})$ be a Poisson manifold and let $B \in \Lambda^2(T ^\ast 
M)$ be the corresponding Poisson tensor. An embedded submanifold $S \subset M $ is
called {\bfi  cosymplectic}
\index{cosymplectic!submanifold}%
\index{submanifold!cosymplectic}%
if 
\begin{description}
\item [(i)] $B^{\sharp}((TS)^{\circ})\cap TS=\{0\}$.
\item [(ii)] $T _sS+ T _s \mathcal{L}_s= T _s M  $,
\end{description}
for any $s \in S $ and $\mathcal{L} _s  $ the symplectic leaf of  $(M,\{ \cdot , \cdot
\})$ containing $s \in S $.
\end{definition}
The cosymplectic submanifolds of a symplectic manifold $(M, \omega)$ are its symplectic
submanifolds. Cosymplectic
submanifolds appear  in the physics literature under the name of {\bfi  second
class constraints}.

\begin{proposition}
\label{properties of cosymplectic submanifolds}
Let $(M, \{ \cdot , \cdot \})$ be a Poisson manifold, $B \in \Lambda^2(T ^\ast 
M)$  the corresponding Poisson tensor, and $S $  a cosymplectic submanifold of $M$.
Then for any $s \in S $
\begin{description}
\item [(i)] $T _s \mathcal{L} _s=(T _s S\cap T _s \mathcal{L} _s)\oplus
B^{\sharp}(s)((T_{s}S)^{\circ}) $, where $\mathcal{L}_s  $ is the symplectic leaf of
$(M,\{ \cdot , \cdot \})$ that contains $s \in S $.
\item [(ii)] $(T _s S)^{\circ}\cap \ker B ^{\sharp} (s)=\{0\}$.
\item [(iii)] $T _s M = B ^{\sharp} (s)((T_{s}S)^{\circ})\oplus T _s S $.
\item  [(iv)]  $B^{\sharp}((TS)^{\circ}) $ is a subbundle of $TM| _S$ and hence $TM|
_S = B^{\sharp}((TS)^{\circ})\oplus TS $.
\item [(v)] The symplectic leaves of $(M,\{ \cdot , \cdot \})$ intersect $S$
transversely and hence $S\cap \mathcal{L} $ is an initial submanifold of $S$, for any
symplectic leaf $\mathcal{L} $ of $(M,\{ \cdot , \cdot \})$.  
\end{description}
\end{proposition}

\begin{theorem}[The Poisson structure of a cosymplectic submanifold]
\label{The Poisson structure of a cosymplectic submanifold theorem}
Let $(M, \{ \cdot , \cdot \})$ be a Poisson manifold, $B \in \Lambda^2(T ^\ast 
M)$  the corresponding Poisson tensor, and $S $  a cosymplectic submanifold of $M$.
Let  $D:= B^{\sharp}((TS)^{\circ}) \subset TM| _S $. Then
\begin{description}
\item [(i)]  $(M, \{ \cdot , \cdot \}, D , S) $ is Poisson reducible.
\item [(ii)] The corresponding quotient manifold equals $S$ and the reduced bracket $\{
\cdot , \cdot \} ^S $ is given by 
\begin{equation}
\label{reduced bracket for cosymplectic manifolds}
\{f,g\} ^S (s)=\{F,G\} (s),
\end{equation}
where $f,g \in C^{\infty}_{S,M}(V) $ are arbitrary and $F,G \in C^{\infty} (U) $ are
local $D$-invariant extensions of $f$ and $g$ around $s\in S$, respectively. 
\item [(iii)] The Hamiltonian vector field $X _f  $ of an arbitrary function $f \in
C^{\infty} _{S,M} (V) $ is given either by
\begin{equation}
\label{expression of the Hamiltonian vector field when invariance}
T i \circ X _f= X _F \circ  i,
\end{equation}
where $F \in C^{\infty} (U) $ is  a local $D$-invariant extension of $f$
and
$i:S\hookrightarrow M $ is the inclusion, or by
\begin{equation}
\label{expression of the Hamiltonian vector field when no invariance}
T i \circ X _f= \pi_S\circ X _{ \overline{F}} \circ  i,
\end{equation}
where $\overline{F} \in C^{\infty} (U) $ is an arbitrary local extension  of $f$ and $\pi_S:TM|
_S \rightarrow TS $ is the projection induced by the Whitney sum decomposition $TM| _S=
B^{\sharp}((TS)^{\circ})\oplus TS $  of  $TM| _S$.
\item [(iv)] The symplectic leaves of $(S,\{ \cdot , \cdot \} ^S)$ are the connected
components of the intersections $S\cap \mathcal{L}$, where $\mathcal{L}$ is
a symplectic leaf of $(M, \{ \cdot , \cdot \})$. Any symplectic leaf of
$(S,\{
\cdot , \cdot \} ^S) $ is a symplectic submanifold of the symplectic leaf of
$(M, \{ \cdot , \cdot \})$ that contains it.
\item [(v)] Let  $\mathcal{L} _s $ and $\mathcal{L} _s^S $  be the symplectic leaves of
$(M, \{ \cdot , \cdot \})$ and $(S,\{ \cdot , \cdot \} ^S)$, respectively, that contain
the point $s \in S $. Let $\omega_{\mathcal{L} _s} $ and $\omega_{\mathcal{L} _s^S} $ be
the corresponding symplectic forms. Then $B^{\sharp}(s)((T_{s}S)^{\circ}) $ is a
symplectic subspace of $T _s \mathcal{L} _s$ and
\begin{equation}
\label{symplectic orthogonal Weinstein}
B^{\sharp}(s)((T_{s}S)^{\circ})= 
\left(T _s \mathcal{L} _s ^S \right)^{\omega_{\mathcal{L}_s} (s)},
\end{equation}
where $\left(T _s \mathcal{L} _s ^S \right)^{\omega_{\mathcal{L}_s} (s)}$
denotes the $\omega_{\mathcal{L}_s} (s)$-orthogonal complement of $T _s
\mathcal{L} _s ^S$ in $T _s \mathcal{L}_s$.

\item [(vi)] Let $B _S \in \Lambda^2(T ^\ast  S)$ be the Poisson tensor associated to
$(S,\{ \cdot , \cdot \} ^S)$. Then
\begin{equation}
\label{something like in Weinstein}
B _S^{\sharp}= \pi_S \circ B ^{\sharp}| _S \circ \pi_S ^\ast, 
\end{equation}
where $\pi_S ^\ast : T ^\ast S \rightarrow  T ^\ast  M| _S $ is the dual of $\pi_S:TM|
_S \rightarrow TS $.
\end{description}
\end{theorem}   
The {\bf Dirac constraints formula} is the expression in coordinates for the bracket
of a cosymplectic submanifold. Let $(M,\{ \cdot , \cdot \})$ be a
$n$-dimensional Poisson manifold and let
$S$ be a $k$-dimensional cosymplectic submanifold of  $M$. Let $z _0 $ be an
arbitrary point in $S$ and $(U, \overline{\kappa})$ a submanifold chart around $z _0 $
such that $\overline{\kappa}=(\overline{\varphi}, \overline{\psi}):U \rightarrow V
_1\times V _2 $, where $V _1 $ and $V _2 $ are two open neighborhoods of the
origin in two Euclidean spaces such that  $\overline{\kappa} (z _0)= \left(
\overline{\varphi}(z _0),
\overline{\psi} (z _0) 
\right)=(0,0) $ and 
\begin{equation}
\label{submanifold propertiy with kappa}
\overline{\kappa}(U\cap S)=V _1\times \{0\}.
\end{equation}
Let $\overline{\varphi}=:(\overline{\varphi} ^1, \ldots,\overline{\varphi} ^k) $ be the
components of $\overline{\varphi} $ and define $\widehat{\varphi} ^1:= \overline{\varphi}
^1|_{U\cap S}, \ldots ,\widehat{\varphi} ^k:= \overline{\varphi}
^k|_{U\cap S}$. Extend  
$\widehat{\varphi} ^1,
\ldots ,\widehat{\varphi} ^k$ to
$D$-invariant functions $ \varphi  ^1 
, \ldots ,\varphi ^k $ on $U$. Since the differentials $\mathbf{d}\widehat{\varphi}
^1(s),
\ldots ,\mathbf{d}\widehat{\varphi} ^k(s)
$ are linearly independent for any $s \in U\cap S $, we can assume (by shrinking $U$ if
necessary) that $\mathbf{d} \varphi 
^1(z),
\ldots ,\mathbf{d} \varphi  ^k(z)
$ are also linearly independent for any $z \in U$. Consequently, $(U, \kappa)$ with
$\kappa:=(\varphi ^1, \ldots, \varphi ^k, \psi ^1, \ldots,\psi^{n-k}) $, is
a submanifold chart for $M$ around $z _0 $ with respect to $S$ such that, by
construction,
\[
\mathbf{d} \varphi^1(s)|_{B^{\sharp}(s)((T_{s}S)^{\circ})}= \cdots=\mathbf{d}
\varphi^k(s)|_{B^{\sharp}(s)((T_{s}S)^{\circ})}=0,
\] 
for any $s \in U\cap S $. This implies that for any $i \in \{1, \ldots, k\} $, $j \in
\{1, \ldots, n-k\} $, and $s \in S $
\[
\{ \varphi^i, \psi ^j\} (s)= \mathbf{d} \varphi ^i (s) \left(X _{\psi ^j}
(s)\right) =0
\]
since $\mathbf{d} \psi ^j (s) \in \left(T _sS \right) ^{\circ} $ by~(\ref{submanifold
propertiy with kappa}) and hence  
\begin{equation}
\label{the x does not belong it 74}
X _{\psi ^j} (s) \in B^{\sharp}(s)((T_{s}S)^{\circ}).
\end{equation}
Additionally, since the functions $\varphi ^1, \ldots, \varphi ^k $ are
$D$-invariant, by~(\ref{expression of the Hamiltonian vector
field when invariance}), it follows that 
\[
X_{\varphi ^1}(s)=X_{\widehat{ \varphi} ^1}(s) \in  T _s S, \ldots, X_{\varphi
^k}(s)=X_{\widehat{ \varphi} ^k}(s) \in  T _s S,
\]
for any $s \in S $. Consequently, $\{X_{\varphi ^1}(s), \ldots, X_{\varphi
^k}(s) , X _{\psi ^1} (s), \ldots ,X _{\psi ^{n-k}} (s)\} $ spans $T _s \mathcal{L} _s $
with 
\[
\{X_{\varphi ^1}(s), \ldots, X_{\varphi
^k}(s)\} \subset T _s S\cap T _s \mathcal{L} _s
\] and 
\[
\{X _{\psi ^1} (s), \ldots ,X
_{\psi ^{n-k}} (s)\} \subset B^{\sharp}(s)((T_{s}S)^{\circ})
.\]
By
Proposition~\ref{properties of cosymplectic submanifolds}{\bf (i)},  
\[
{\rm
span}\{X_{\varphi ^1}(s), \ldots, X_{\varphi ^k}(s)\}=T _s S\cap T _s \mathcal{L} _s
\]
and 
\[
{\rm span}\{X _{\psi ^1} (s), \ldots ,X _{\psi ^{n-k}} (s)\} =
B^{\sharp}(s)((T_{s}S)^{\circ}). 
\]
Since 
$
\dim \left( B^{\sharp}(s)((T_{s}S)^{\circ})\right)=n-k
$
by Proposition~\ref{properties of cosymplectic
submanifolds}{\bf (iii)},  it follows that $\{X _{\psi ^1} (s), \ldots ,X
_{\psi ^{n-k}} (s)\}$ is a basis of $B^{\sharp}(s)((T_{s}S)^{\circ}) $.

Since $B^{\sharp}(s)((T_{s}S)^{\circ})$ is a symplectic subspace of $T_s
\mathcal{L}_s$ by Theorem~\ref{The Poisson structure of a
cosymplectic submanifold theorem}{\bf (v)}, there exists some $r \in
\mathbb{N} $ such that  $n-k=2r
$ and, additionally, the matrix $C (s) $ with entries 
\[
C^{ij}(s):=\{ \psi^i, \psi^j\} (s), \qquad i,j \in \{1, \ldots, n-k\}
\]
is invertible. Therefore, in the coordinates $(\varphi ^1, \ldots, \varphi^k, \psi^1,
\ldots, \psi^{n-k}) $, the matrix associated to the Poisson tensor $B (s) $ is
\[
B (s)=\left(
\begin{array}{cc}
B _S (s)& 0\\
0 &C (s)
\end{array}
\right),
\]
\normalsize
where $B _S \in \Lambda^2(T ^\ast  S)$ is the Poisson tensor associated to
$(S,\{ \cdot , \cdot \} ^S)$. Let $C_{ij}(s)$ be the entries of the matrix $C(s)^{-1}$.

\begin{proposition}[Dirac formulas]
\index{Dirac!formula}%
\index{formula!Dirac}%
In the coordinate neighborhood $(\varphi ^1, \ldots, \varphi^k, \psi^1,
\ldots, \psi^{n-k}) $ constructed above and for $s \in S $ we have, for any $f,g \in 
C^{\infty}_{S,M} (V) $:
\begin{equation}
\label{dirac formula vector field 351}
X _f (s)=X _F (s)-\sum_{i,j=1}^{n-k}\{F, \psi ^i\} (s) C_{ij}(s) X_{\psi ^j}(s)
\end{equation}
and
\begin{equation}
\label{dirac formula bracket 351}
\{f,g\} ^S (s)=\{F,G\}(s)-\sum_{i,j=1}^{n-k}\{F, \psi ^i\} (s) C_{ij}(s) \{\psi
^j,G\}(s),
\end{equation}
where $F,G \in C^{\infty} (U) $ are arbitrary local extensions of $f$ and $ g $,
respectively, around $s \in S $.
\end{proposition}

\end{document}